\documentclass[a4paper,12pt, oneside, reqno]{amsart}
\usepackage{mathrsfs}
\usepackage{amsfonts}
\usepackage{amsmath,amstext,amssymb,amsthm}
\usepackage[polish, english]{babel}
\usepackage{amsxtra}
\usepackage{color}
\usepackage{dsfont}
\usepackage[margin=2cm, centering]{geometry}
\usepackage[bbgreekl]{mathbbol}
\usepackage{soul}
\usepackage{graphics}

\usepackage[colorlinks,citecolor=blue,urlcolor=blue,bookmarks=true]{hyperref}

\DeclareSymbolFontAlphabet{\mathbb}{AMSb}
\DeclareSymbolFontAlphabet{\mathbbl}{bbold}

\usepackage{accents}

\theoremstyle{plain}
\newtheorem{theorem}{Theorem}[section] 
\newtheorem{lemma}[theorem]{Lemma}
\newtheorem{proposition}[theorem]{Proposition}
\newtheorem{corollary}[theorem]{Corollary}

 \theoremstyle{definition}
\newtheorem{defn}{Definition}

\newtheorem*{remark*}{Remark} 
\newtheorem{remark}[theorem]{Remark}

\newcommand{\R}{\mathbb{R}}
\newcommand{\Rd}{{\R^{d}}}

\newcommand{\ind}{\mathds{1}}
\renewcommand{\leq}{\leqslant}
\renewcommand{\geq}{\geqslant}

\def\ind{{\bf 1}}
\def\qed{{\hfill $\Box$ \bigskip}}
\def\LL{{\mathcal L}}

\def\N{{\mathbb N}}
\def\pf{\noindent{\bf Proof.} }

\def\({\left(} 
\def\){\right)} 
\def\[{\left[}
\def\]{\right]} 
\def\<{\langle} 
\def\>{\rangle}

\def \Pa{\rm{(P1)}}
\def \Pb{\rm{(P2)}}
\def \Pc{\rm{(P3)}}

\def \Qa{\rm{(Q1)}}
\def \Qb{\rm{(Q2)}}

\def \GC_r{\rm{(PQ)}}

\newcommand{\lah}{\alpha_h}
\newcommand{\uah}{\beta_h}

\newcommand{\err}[2]{\rho_{#1}^{#2}} 
\newcommand{\param}{\sigma}  
\newcommand{\parame}{\sigma_e}  
\newcommand{\lmCJ}{\gamma_0}  

\newcommand{\aF}{\mathcal{F}} 

\newcommand{\rr}{\Upsilon} 

\newcommand{\ka}{\kappa_3}  
\newcommand{\kb}{\kappa_4}  

\definecolor{ks}{rgb}{0.7,0.1,0.2}

\allowdisplaybreaks

\title{Regularity of fundamental solutions for L{\'e}vy-type operators}

\thanks{The research was partially supported by
National Science Centre (Poland)
grant 2016/23/B/ST1/01665.}

\author[K. Szczypkowski]{Karol Szczypkowski}
\address{
	Karol Szczypkowski\\
	Wydzia{\l} Matematyki,
	Politechnika Wroc{\l}awska\\
	Wyb. Wyspia\'{n}skiego 27\\
	50-370 Wroc{\l}aw\\
	Poland
}
\email{karol.szczypkowski@pwr.edu.pl}

\date{}

\begin{document}

\begin{abstract}
For a class of non-symmetric non-local L{\'e}vy-type operators
$\LL^{\kappa}$, which include those of the form
$$
\LL^{\kappa}f(x):= \int_{\R^d}( f(x+z)-f(x)- \ind_{|z|<1} \left<z,\nabla f(x)\right>)\kappa(x,z)J(z)\, dz\,,
$$
we prove 
regularity of the fundamental solution $p^{\kappa}$ to the equation  $\partial_t =\LL^{\kappa}$.
\end{abstract}

\maketitle



\section{Introduction}

L{\'e}vy-type processes
are stochastic models
that can be used to approximate physical, biological or financial phenomena.
A local  expansion of such a process is described by
its infinitesimal generator which is a L{\'e}vy-type operator $\LL$. 
The non-local part of that operator is responsible for and describes the jumps of the process.
In the recent years  a lot of effort has been put into understanding
purely non-local L{\'e}vy-type operators.
At first the case with constant coefficients attracted most of the attention,
but the literature concerning non-constant coefficient is growing rapidly,
including 
\cite{MR3652202},
\cite{MR3500272}, 
\cite{MR3817130},
\cite{GS-2018}
\cite{FK-2017},
\cite{MR3765882},
\cite{BKS-2017}, 
\cite{KR-2017},
\cite{PJ},
\cite{S-2018},
\cite{MR3294616},
\cite{MR2163294}, 
\cite{MR3353627},
\cite{MR3500272}.
The parametrix method \cite{zbMATH02644101}, \cite{MR2093219} used in those papers leads to a construction of
the fundamental solution of the equation $\partial_t u = \LL u$
or the heat kernel of the process that is a unique solution
to the martingale problem for $\LL$.

The subject of the present paper
is  non-local L{\'e}vy-type operators with H{\"o}lder continuous coefficients 
considered in \cite{GS-2018} and \cite{S-2018}
(see Definition~\ref{def-r:GC_r}). 
A typical example here is the operator
\begin{align*}
\LL f(x)= \int_{\Rd} \left( f(z+x)-f(x) -  
\ind_{|z|<1} \left<z,\nabla f(x)\right>\right)\frac{\kappa(x,z)}{|z|^{d+\alpha}}\, dz\,,
\end{align*}
where $\alpha\in (0,2)$, $\kappa$ is bounded from below and above by positive constants, and $\beta$-H{\"o}lder continuous in the first variable with $\beta \in (0,1)$.
The usual result concerning the regularity of the fundamental solution
of $\partial_t = \LL $
is $\gamma$-H{\"o}lder continuity with $\gamma<\alpha$.
We improve that result in more general setting
taking into account the $\beta$ regularity  of the coefficient $\kappa$.

Of particular interest are the existence, estimates and regularity of the gradient of the fundamental solution for L{\'e}vy or L{\'e}vy-type operators \cite{MR2995789},\cite{MR3413864}, \cite{MR3981133}, \cite{GS-2017}; \cite{MR3544166}, \cite{MR3472835}. In this context our assumptions will imply $\alpha>1/2$. In a recent paper
\cite{LSX-2020} this restriction was removed at the expense of additional constraints on the coefficient $\kappa$ requiring strong symmetry properties in~$z$. 
As already mentioned, the purpose of the present paper is to cover a wide class of operators and coefficients discussed in \cite{GS-2018} and \cite{S-2018}. In particular ones that are not symmetric in the $z$ variable, thus not considered in \cite{LSX-2020}.
What is more, such non-symmetry may cause a (time dependent) non-zero {\it internal drift} with a coefficient that is unbounded as time tends to zero (see \cite[Example~2]{S-2018}), which we cover in case $\Qa$ below.

Under certain conditions, which require $\alpha>1$, we also prove existence, estimates and regularity of the second derivatives of the fundamental solution. To the best of the author's knowledge, this result is new  even for the operator $\LL$ above 
 under any assumptions on the coefficient $\kappa$
 that is not constant in $x$.

\section{The setting and main results}

Let $d\in\N$ and
$\nu:[0,\infty)\to[0,\infty]$ be a non-increasing  function  satisfying
$$\int_{\Rd}  (1\land |x|^2) \nu(|x|)dx<\infty\,.$$
For real numbers $a$ and $b$ we write as usual $a \land b = \min (a,b)$ and $a\vee b = \max(a,b)$.
Consider
$J: \Rd  \to [0, \infty]$ 
 such that for some  $\lmCJ \in [1,\infty)$ and 
 all $x\in \Rd$,
\begin{equation}\label{e:psi1}
\lmCJ^{-1} \nu(|x|)\leq J(x) \leq \lmCJ \nu(|x|)\,.
\end{equation}
Further, suppose that 
$\kappa(x,z)$ is a Borel 
function on $\R^d\times \Rd$ such that
\begin{equation}\label{e:intro-kappa}
0<\kappa_0\leq \kappa(x,z)\leq \kappa_1\, , 
\end{equation}
and 
for some $\beta\in (0,1)$,
\begin{equation}\label{e:intro-kappa-holder}
|\kappa(x,z)-\kappa(y,z)|\leq \kappa_2|x-y|^{\beta}\, .
\end{equation}
For $r>0$ we define 
$$
h(r):= \int_0^\infty \left(1\land \frac{|x|^2}{r^2}\right) \nu(|x|)dx\,,\qquad \quad
K(r):=r^{-2} \int_{|x|<r}|x|^2 \nu(|x|)dx\,.
$$
The above functions  play a prominent role in the paper.
Our main assumption is \emph{the weak  scaling condition} at the origin:  there exist $\lah\in (0,2]$ and $C_h \in [1,\infty)$ such that 
\begin{equation}\label{eq:intro:wlsc}
h(r)\leq C_h\,\lambda^{\lah}\,h(\lambda r)\, ,\quad \lambda\leq 1, r\leq 1\,,
\end{equation}
and in a  similar fashion, there exist $\uah\in (0,2]$ and $c_h\in (0,1]$ such that
\begin{equation}\label{eq:intro:wusc}
 h(r)\geq c_h\,\lambda^{\uah}\,h(\lambda r)\, ,\quad \lambda\leq 1, r\leq 1\, .\\
\end{equation}
Furthermore,
suppose there are (finite) constants $\ka, \kb\geq 0$ such that
\begin{align}\label{e:intro-kappa-crit}
 \sup_{x\in\Rd}\left| \int_{r\leq |z|<1} z\, \kappa(x,z) J(z)dz \right| &\leq \ka  rh(r)\,,
 \qquad r\in (0,1],\\
\left| \int_{r\leq |z|<1} z\, \big[ \kappa(x,z)- \kappa(y,z)\big] J(z)dz \right| &\leq \kb |x-y|^{\beta} rh(r)\,, \qquad r\in (0,1]. \label{e:intro-kappa-crit-H}
\end{align}

\begin{defn}\label{def-r:GC_r}
We say that $\GC_r$ holds if one of the following sets of assumptions is satisfied,
\begin{enumerate}
\item[] 
\begin{enumerate}
\item[$\Pa$] \quad  \eqref{e:psi1}--\eqref{eq:intro:wlsc}
hold and  $1< \lah \leq 2$;
\item[$\Pb$] \quad \eqref{e:psi1}--\eqref{eq:intro:wusc}
hold and   $0<\lah \leq \uah <1$;
\item[$\Pc$] \quad \eqref{e:psi1}--\eqref{eq:intro:wlsc} hold, $J$ is symmetric and $\kappa(x,z)=\kappa(x,-z)$, $x,z\in\Rd$;
\item[] 
\item[$\Qa$] \quad \eqref{e:psi1}--\eqref{eq:intro:wlsc}
hold, $\lah=1$; \eqref{e:intro-kappa-crit} and \eqref{e:intro-kappa-crit-H} hold;
\item[$\Qb$] \quad \eqref{e:psi1}--\eqref{eq:intro:wusc}
hold,  $0<\lah \leq \uah <1$ and $1-\lah<\beta \land \lah$; \eqref{e:intro-kappa-crit} and \eqref{e:intro-kappa-crit-H} hold.
\end{enumerate}
\end{enumerate}
\end{defn}

Our aim is to prove regularity of the heat kernel $p^{\kappa}$
of a non-local non-symmetric L{\'e}vy-type operator~$\LL^{\kappa}$,
i.e., of a fundamental solution to the equation $\partial_t u =\LL^{\kappa} u$.
For each of cases $\Pa$, $\Qa$, $\Qb$ the operator
under consideration is defined as
\begin{align}
\LL^{\kappa}f(x)&:= \int_{\Rd}( f(x+z)-f(x)- \ind_{|z|<1} \left<z,\nabla f(x)\right>)\kappa(x,z)J(z)\, dz \,. \label{e:intro-operator-a1}
\end{align}
If $\Pb$ holds we consider
\begin{align}
\LL^{\kappa}f(x)&:= \int_{\Rd}( f(x+z)-f(x))\kappa(x,z)J(z)\, dz\,. \label{e:intro-operator-a2}
\end{align}
If $\Pc$ holds we discuss
\begin{align}
\LL^{\kappa}f(x)&:= \frac1{2}\int_{\Rd}( f(x+z)+f(x-z)-2f(x))\kappa(x,z)J(z)\, dz \,.\label{e:intro-operator-a3}
\end{align}

It was shown in \cite[Theorem~1.1]{GS-2018} and \cite[Theorem~2.1]{S-2018}
 that under $\GC_r$ the function  $p^{\kappa}$  exists and is unique within a certain class of functions.
In fact, the heat kernel was constructed using the Levi's paramterix method, i.e.,
\begin{align}\label{e:p-kappa}
p^{\kappa}(t,x,y)=
p^{\mathfrak{K}_y}(t,x,y)+\int_0^t \int_{\Rd}p^{\mathfrak{K}_z}(t-s,x,z)q(s,z,y)\, dzds\,,
\end{align}
where $q(t,x,y)$ solves the equation
\begin{align*}
q(t,x,y)=q_0(t,x,y)+\int_0^t \int_{\Rd}q_0(t-s,x,z)q(s,z,y)\, dzds\,,
\end{align*}
and $q_0(t,x,y)=\big(\LL_x^{{\mathfrak K}_x}-\LL_x^{{\mathfrak K}_y}\big) p^{\mathfrak{K}_y}(t,x,y)$.
The function $p^{\mathfrak{K}_w}$ is the heat kernel of the L{\'e}vy operator $\LL^{\mathfrak{K}_w}$
obtained from the operator $\LL^{\kappa}$ by freezing its coefficient: $\mathfrak{K}_w(z)=\kappa(w,z)$.
For $t>0$ and $x\in \R^d$ we define {\it the bound function}
\begin{equation}\label{e:intro-rho-def}
\rr_t(x):=\left( [h^{-1}(1/t)]^{-d}\land \frac{tK(|x|)}{|x|^{d}} \right) .
\end{equation}
We refer the reader to
\cite[Theorems~1.2 and~1.4]{GS-2018} and
\cite[Theorems~2.2 and~2.4]{S-2018} 
for a collection of properties of the function $p^{\kappa}$, including estimates, H{\"o}lder continuity, differentiability and gradient estimates.
For instance,
under $\GC_r$
for all $T>0$ and $\gamma \in [0,1] \cap[0,\lah)$,
there is  $c>0$ such that for all $t\in (0,T]$ and $x,x',y\in \Rd$,
\begin{align*}
\left|p^{\kappa}(t,x,y)-p^{\kappa}(t,x',y)\right| \leq c 
 (|x-x'|^{\gamma}\land 1) \left[h^{-1}(1/t)\right]^{-\gamma} \big( \rr_t(y-x)+ \rr_t(y-x') \big).
\end{align*}
Recall also that \cite[Theorem~1.2 (6)]{GS-2018} and
\cite[Theorem~2.2 (6)]{S-2018} provide an upper bound for $|\nabla_{x} p^{\kappa}(t,x,y)|$ under conditions $\GC_r$ and $\lah+\beta\land \lah>1$.

\vspace{\baselineskip}

Here are the main results of the present paper.
For the meaning of $\parame$ see Definition~\ref{def:parame}.

\begin{theorem}\label{thm-r:1}
Assume $\GC_r$.
Let 
$r_0 \in [0,1]\cap [0,\lah+\beta\land \lah)$.
For every $T>0$ there exists a constant
$c=c(d,T,\parame,r_0)$ such that for all $t\in (0,T]$,
$x,x',y\in\Rd$ and $r\in [0,r_0]$,
\begin{align*}
|p^{\kappa}(t,x,y)-p^{\kappa}(t,x',y)| \leq c \,(|x-x'|^{r}\land 1) \left[ h^{-1}(1/t)\right]^{-r}
 \big( \rr_t(y-x') + \rr_t(y-x)\big)\,.
\end{align*}
\end{theorem}

\begin{theorem}\label{thm-r:2a}
Assume $\GC_r$ and suppose that
$\lah+\beta\land \lah>1$. 
Let $r_0 \in [0,1]\cap [0,\lah+\beta\land \lah-1)$.
For every $T>0$ there exists a constant
$c=c(d,T,\parame,r_0)$ such that for all $t\in (0,T]$,
$x,x',y\in\Rd$ and $r\in [0,r_0]$,
\begin{align*}
|\nabla_{x} p^{\kappa}(t,x,y)-\nabla_{x'} p^{\kappa}(t,x',y)| \leq c \,(|x-x'|^{r}\land 1) \left[ h^{-1}(1/t)\right]^{-1-r}
 \big( \rr_t(y-x') + \rr_t(y-x)\big)\,.
\end{align*}
\end{theorem}

\begin{theorem}\label{thm-r:2}
Assume $\GC_r$ and suppose that
$\lah+\beta\land \lah>2$.
Let $r_0 \in [0,1]\cap [0,\lah+\beta\land \lah-2)$.
For every $T>0$ there exists a constant
$c=c(d,T,\parame,r_0)$ such that for all $t\in (0,T]$,
$x,x',y\in\Rd$ and $r\in [0,r_0]$,
\begin{align*}
|\nabla_x^2\, p^{\kappa}(t,x,y)|\leq c  \left[h^{-1}(1/t)\right]^{-2} \rr_t(y-x)\,,
\end{align*}
\begin{align*}
|\nabla_{x}^2 p^{\kappa}(t,x,y)-\nabla_{x'}^2 p^{\kappa}(t,x',y)| \leq c \,(|x-x'|^{r}\land 1) \left[ h^{-1}(1/t)\right]^{-2-r}
 \big( \rr_t(y-x') + \rr_t(y-x)\big)\,.
\end{align*}
\end{theorem}

\vspace{\baselineskip}
The condition $\lah+\beta\land \lah>1$ equivalently means that there is $\beta_1 \in [0,\beta]\cap [0,\lah)$ such that
$\lah+\beta_1>1$, and it may hold only if $\lah>1/2$.
Similarly, $\lah+\beta\land \lah>2$
is equivalent to the existence of $\beta_1 \in [0,\beta]\cap [0,\lah)$
such that $\lah+\beta_1>2$, and it requires $\lah>1$. 
We  note that even the existence of second derivatives of $p^{\kappa}$
in Theorem~\ref{thm-r:2} is a new result.

\begin{defn}\label{def:parame}
Following \cite{GS-2018},
in the case
$\Pa$, $\Pb$, $\Pc$ we respectively consider the set of parameters 
$\param_1 = (\lmCJ,\kappa_0,\kappa_1,\lah, C_h,h)$,
$\param_2= (\lmCJ,\kappa_0,\kappa_1,\lah,\uah, C_h, c_h,h)$,
$\param_3= (\lmCJ,\kappa_0,\kappa_1,\lah,C_h,h)$,
which we abbreviate to $\param$.
Similarly, after \cite{S-2018}
we put $\param=(\lmCJ,\kappa_0,\kappa_1,\ka,\lah, C_h,h)$
under $\Qa$ or $\Qb$.
We extend the set of parameters $\param$ to $\parame$
by adding constant $\kappa_2$ in the cases $\Pa$, $\Pb$, $\Pc$,
and $\kappa_2,\kb$ in the cases $\Qa$, $\Qb$.
Abusing the notation we have
$\parame=(\sigma,\kappa_2)$ under $\Pa$, $\Pb$, $\Pc$,
and $\parame=(\sigma,\kappa_2,\kb)$ under $\Qa$, $\Qb$.
\end{defn}

In the whole paper {\bf we assume that $\GC_r$ holds.}

\section{Preliminaries}

To prove the main results we will obviously use the representation \eqref{e:p-kappa} of $p^\kappa(t,x,y)$. This formula consists of  $p^{\mathfrak{K}_y}(t,x,y)$, which is known in the literature as {\it the zero order approximation}, and the integral part called {\it the remainder}.
In the whole paper we follow and use consistently the notation of \cite{GS-2018} and \cite{S-2018}.
To study the remainder we let
\begin{equation}\label{e:phi-y-def}
\phi_y(t,x,s):=\int_{\Rd} p^{\mathfrak{K}_z}(t-s,x,z)q(s,z,y)\, dz, \quad x \in \Rd, \, 
0< s<t\,,
\end{equation}
and
\begin{equation}\label{e:def-phi-y-2}
\phi_y(t,x):=\int_0^t \phi_y(t,x,s)\, ds =\int_0^t \int_{\Rd}p^{\mathfrak{K}_z}(t-s,x,z)q(s,z,y)\, dzds\,.
\end{equation}
See  \cite[Theorem~3.7]{GS-2018}
and \cite[Theorem~6.2]{S-2018}
for the definition of $q(t,x,y)$ as a series. In our proofs we only use the properties of $q(t,x,y)$ (mostly already known in \cite{GS-2018} and \cite{S-2018}) rather than its concrete structure.
We start by investigating the regularity of the zero order approximation and we do so in a slightly more general context that will be useful when dealing with the remainder.
We introduce the following expression: for $t>0$ and $x,y,z\in\Rd$,
\begin{align*}
\aF_{2}(t,x,y;z)&:=\rr_t(y-x-z)\ind_{|z|\geq h^{-1}(1/t)}+ \left[ \left(\frac{|z|}{h^{-1}(1/t)}\right)\wedge 1\right] \rr_t(y-x).
\end{align*}
In \cite{GS-2018} and \cite{S-2018} also $\aF_{1}$ was considered, but we will not need it in our analysis.
In what follows functions, $\mathfrak{K}$
and $p^{\mathfrak{K}}$
are as in \cite[Section~2]{GS-2018} and \cite[Section~4]{S-2018}.
Namely,
$\mathfrak{K}\colon \Rd \to [0,\infty)$ is such that
$$
0<\kappa_0 \leq \mathfrak{K}(z) \leq \kappa_1\,.
$$
Now, we consider an operator 
$\LL^{\mathfrak{K}}$
defined by taking $\kappa(x,z)=\mathfrak{K}(z)$ in
\eqref{e:intro-operator-a1} for $\Pa$, $\Qa$, $\Qb$;  
\eqref{e:intro-operator-a2} for $\Pb$; and  \eqref{e:intro-operator-a3} for $\Pc$. The operator uniquely determines a L{\'e}vy process and its density $p^{\mathfrak{K}}(t,x,y)=p^{\mathfrak{K}}(t,y-x)$, see \cite[Section~6]{GS-2018}.

\begin{lemma}\label{lem-r:increments}
For every $T>0$ there exists a constant
$c=c(d,T,\param)$
such that for all $t\in(0,T]$ and $x,y,z \in \Rd$,
\begin{align}\label{ineq-r:est_diff_1}
|p^{\mathfrak{K}}(t,x+z,y)-p^{\mathfrak{K}}(t,x,y) | 
&\leq c\, \aF_{2}(t,x,y;z),\\
\label{ineq-r:est_diff_grad_1}
|\nabla_{x} p^{\mathfrak{K}}(t,x+z,y)- \nabla_x p^{\mathfrak{K}}(t,x,y) | 
&\leq c \left[h^{-1}(1/t)\right]^{-1} \aF_{2}(t,x,y;z),\\
\label{ineq-r:est_diff_grad_2}
|\nabla_{x}^2\, p^{\mathfrak{K}}(t,x+z,y)- \nabla_x^2\, p^{\mathfrak{K}}(t,x,y) | 
&\leq c \left[h^{-1}(1/t)\right]^{-2} \aF_{2}(t,x,y;z).
\end{align}
\end{lemma}
\pf
The inequalities follow from
\cite[Proposition~2.1]{GS-2018} and
\cite[Proposition~4.1]{S-2018}
that provide upper bounds for derivatives of $p^{\mathfrak{K}}$ in spatial variable. 
If $|z|\geq h^{-1}(1/t)$ we bound each term of the difference separately using these bounds. If $|z|< h^{-1}(1/t)$ 
it suffices to write the increment as an integral of the   derivative and then use the bounds while removing small shifts in arguments using \cite[Corollary~5.10]{GS-2018};  cf. the proofs of \cite[Lemma~2.3]{GS-2018}, \cite[Lemma~4.3]{S-2018}.
\qed

Inequalities \eqref{ineq-r:est_diff_1}, \eqref{ineq-r:est_diff_grad_1}, \eqref{ineq-r:est_diff_grad_2} can be  written equivalently as
\begin{align*}
\left|p^{\mathfrak{K}}(t,x',y)-p^{\mathfrak{K}}(t,x,y)\right|&\leq c \left(\frac{|x'-x|}{h^{-1}(1/t)} \land 1\right) \big( \rr_t(y-x') + \rr_t(y-x)\big)\,, \\
\left|\nabla_x p^{\mathfrak{K}}(t,x',y)-\nabla_x p^{\mathfrak{K}}(t,x,y)\right|&\leq c
\left(\frac{|x'-x|}{h^{-1}(1/t)} \land 1\right) \left[h^{-1}(1/t)\right]^{-1}\big( \rr_t(y-x') + \rr_t(y-x)\big)\,,\\
\left|\nabla_x^2\, p^{\mathfrak{K}}(t,x',y)-\nabla_x^2\, p^{\mathfrak{K}}(t,x,y)\right|&\leq c
\left(\frac{|x'-x|}{h^{-1}(1/t)} \land 1\right) \left[h^{-1}(1/t)\right]^{-2}\big( \rr_t(y-x') + \rr_t(y-x)\big)\,.
\end{align*}

\begin{corollary}\label{cor-r:est_diff_grad}
For every $T>0$ there exists a constant
$c=c(d,T,\param)$
such that for all $t\in(0,T]$, $x,x',y \in \Rd$ and $\gamma\in [0,1]$,
\begin{align*}
|p^{\mathfrak{K}}(t,x',y)-p^{\mathfrak{K}}(t,x,y) | 
\leq c (|x-x'|^{\gamma}\land 1) \left[h^{-1}(1/t)\right]^{-\gamma} 
 \big( \rr_t(y-x') + \rr_t(y-x)\big),
\end{align*}
\begin{align*}
|\nabla_{x'} p^{\mathfrak{K}}(t,x',y)- \nabla_{x} p^{\mathfrak{K}}(t,x,y) | 
\leq c (|x-x'|^{\gamma}\land 1) \left[h^{-1}(1/t)\right]^{-\gamma-1} 
 \big( \rr_t(y-x') + \rr_t(y-x)\big),
\end{align*}
\begin{align*}
|\nabla_{x'}^2\, p^{\mathfrak{K}}(t,x',y)- \nabla_{x}^2\, p^{\mathfrak{K}}(t,x,y) | 
\leq c (|x-x'|^{\gamma}\land 1) \left[h^{-1}(1/t)\right]^{-\gamma-2} 
 \big( \rr_t(y-x') + \rr_t(y-x)\big).
\end{align*}
\end{corollary}

\vspace{\baselineskip}
Corollary~\ref{cor-r:est_diff_grad} already covers the targeted regularity for the zero order approximation, because it can be applied to $\mathfrak{K}_w(z)=\kappa(w,z)$ to get inequalities which are uniform in $w\in\Rd$, hence in particular uniform for the family of functions $\{p^{\mathfrak{K}_y}(t,x,y)\colon  t\in (0,T],\, x,y\in\Rd \}$ .

Now, we focus on the remainder \eqref{e:def-phi-y-2}. We will clearly have to use estimates for $q(t,x,y)$ and 
again  inequalities \eqref{ineq-r:est_diff_1}--\eqref{ineq-r:est_diff_grad_2}. However, if we simply  apply these two properties, we immediately run into the  problem of integrability of $[h^{-1}(1/(t-s))]^{-\tilde{r}}$ with respect to $s\in [t/2,t)$, imposing constraints on the value of $\tilde{r}$, which does not lead to our main results. Therefore, we have to find and exploit cancellations that take place in the integrals defining \eqref{e:def-phi-y-2}. 
We write
\begin{align*}
\phi_y(t,x)-\phi_y(t,x')
= \int_0^{t/2} \left( \phi_y(t,x,s)-\phi_y(t,x',s)\right) ds
+ \int_{t/2}^t \left( \phi_y(t,x,s)-\phi_y(t,x',s)\right) ds
\end{align*}
and we treat the second term as follows: 
\begin{align*}
\int_{t/2}^t ( \phi_y(t,x,s) &-  \phi_y(t,x',s)) ds
=\int_{t/2}^t \int_{\Rd} \left( p^{\mathfrak{K}_z}(t-s,x,z)- p^{\mathfrak{K}_z}(t-s,x',z)\right) q(s,z,y)\, dzds \\
&= \int_{t/2}^t \int_{\Rd} \left( p^{\mathfrak{K}_z}(t-s,x,z)-p^{\mathfrak{K}_z}(t-s,x',z) \right) \big( q(s,z,y)-q(s,x,y) \big) \,dzds \\
&\quad +\int_{t/2}^t\int_{\Rd} \left( p^{\mathfrak{K}_z}(t-s,x,z)-p^{\mathfrak{K}_z}(t-s,x',z) \right) dz \,\,  q(s,x,y) \,ds\,.
\end{align*}
Similar, though slightly different, decompositions are used for the first and second order derivatives under appropriate assumptions.
Roughly speaking, the goal is achieved by making use of the regularity of the coefficient $\kappa(x,z)$ in $x$, and the so called convolution inequalities that involve space or space-time integrals of the following functions
for certain $\gamma,\beta\in \R$:
$$
\err{\gamma}{\beta}(t,x):= \left[h^{-1}(1/t)\right]^{\gamma} \left(|x|^{\beta}\land 1\right) t^{-1} \rr_t(x)\,.
$$
The inequalities are collected in \cite[Lemma~5.17]{GS-2018}, and we will simply refer to that result  whenever using them. 
The very initial step to detect cancellations uses the inequalities in \cite[Theorem~2.11]{GS-2018} and \cite[Proposition~5.3]{S-2018}.
We also note that for the second derivatives we have to prove more, because this case was less studied in \cite{GS-2018} or \cite{S-2018}

\begin{remark}\label{rem-r:monot_h}
We often use the monotonicity of the function $h$ and $h^{-1}$, see \cite[Lemma~5.1]{GS-2018}, in particular $[h^{-1}(1/(t-s))]^{\gamma} \leq [h^{-1}(1/t)]^{\gamma}\leq h^{-1}(1/T)\vee 1$ holds for all $0<s<t\leq T$, $\gamma\in [0,1]$.
\end{remark}

\begin{remark}\label{rem-r:stretch}
Certain technical results of \cite{GS-2018}, e.g., 
\cite[Lemma~5.3]{GS-2018} with $u=1/t$,
in view of~\eqref{eq:intro:wlsc}
provide inequalities that hold for $t<1/h(1)$.
Using \cite[Remark~5.2]{GS-2018}
we may extend those inequalities to hold for $t\in (0,T]$
by increasing the constant $C_h$ to $C_h [h^{-1}(1/T)\vee 1]^2$.
\end{remark}

We will also need a slight improvement of \cite[(38)]{GS-2018} and \cite[(54)]{S-2018} concerning
the dependence of the constant on the parameter $\gamma>0$.

\begin{lemma}\label{lem-r:q_reg_impr}
Let $\beta_1\in (0,\beta]\cap (0,\lah)$ and $0<\gamma_1\leq \beta_1$.
For every $T>0$ there exists a constant
$c=c(d,T,\parame,\beta_1,\gamma_1)$
such that for all $t\in(0,T]$, $x,x',y \in \Rd$ and $\gamma\in [\gamma_1,\beta_1]$,
\begin{align*}
&|q(t,x,y)-q(t,x',y)|\nonumber\\
&\leq c \left(|x-x'|^{\beta_1-\gamma}\land 1\right)
\left\{\big(\err{\gamma}{0}+\err{\gamma-\beta_1}{\beta_1}\big)(t,x-y)+\big(\err{\gamma}{0}+\err{\gamma-\beta_1}{\beta_1}\big)(t,x'-y)\right\}\,.
\end{align*}
\end{lemma}
\pf
This formulation of the H{\"o}lder continuity of $q$ has the same proof as
\cite[(38)]{GS-2018}.
One only needs to pay attention to explicit constants when applying
\cite[Lemma~5.17(c)]{{GS-2018}}. In particular, the monotonicity of the Beta function is used. See also Remark~\ref{rem-r:monot_h}.
\qed

\section{Regularity - part I}\label{sec:part1}

We start with several technical lemmas before we prove the key Proposition~\ref{prop-r:key_0}.
\begin{lemma}\label{lem-r:for_cancel_0}
For every $T>0$ there exists a constant
$c=c(d,T,\parame)$
such that for all $t\in(0,T]$ and $x,x',y,w,w' \in \Rd$
satisfying $|x-x'|\leq h^{-1}(1/t)$,
\begin{align*}
\left| p^{\mathfrak{K}_w}(t,x,y)- p^{\mathfrak{K}_w}(t,x',y)
-\left( p^{\mathfrak{K}_{w'}}(t,x,y)- p^{\mathfrak{K}_{w'}}(t,x',y)\right)\right|&\\
 \leq 
c \left( \frac{|x-x'|}{h^{-1}(1/t)}\right) (|w-w'|^{\beta}\land 1) & \,\rr_{t}(y-x)\,.
\end{align*}
\end{lemma}

\pf
Let $w_0\in\Rd$  be fixed. Define $\mathfrak{K}_0(z)=(\kappa_0/(2\kappa_1)) \kappa(w_0,z)$ and $\widehat{\mathfrak{K}}_w (z)=\mathfrak{K}_w(z)- \mathfrak{K}_0(z)$. By the construction of the L{\'e}vy process we have
\begin{align}\label{eq-r:przez_k_0-impr}
p^{\mathfrak{K}_w}(t,x,y)=\int_{\Rd} p^{\mathfrak{K}_0}(t,x,\xi)
p^{\widehat{\mathfrak{K}}_w}(t,\xi,y)\,d\xi\,.
\end{align}
Thus,
\begin{align*}
&p^{\mathfrak{K}_w}(t,x,y) - p^{\mathfrak{K}_w}(t,x',y)
-\left(p^{\mathfrak{K}_{w'}}(t,x,y)-p^{\mathfrak{K}_{w'}}(t,x',y)\right)\\
&\quad = \int_{\Rd} \left( p^{\mathfrak{K}_0}(t, x,\xi)- p^{\mathfrak{K}_0}(t, x',\xi) \right)  \left( p^{\widehat{\mathfrak{K}}_w}(t,\xi,y)-p^{\widehat{{\mathfrak{K}}}_{w'}}(t,\xi,y)\right) d\xi\,.
\end{align*}
By \eqref{ineq-r:est_diff_1}
and
\cite[Theorem~2.11]{GS-2018}, \cite[Proposition~5.3]{S-2018},
for $|x-x'|\leq h^{-1}(1/t)$ we get
$$
\left| p^{\mathfrak{K}_0}(t, x,\xi)- p^{\mathfrak{K}_0}(t, x',\xi) \right| \leq c \left( \frac{|x-x'|}{h^{-1}(1/t)}\right) \rr_t(\xi-x)\,,
$$
$$
\left| p^{\widehat{\mathfrak{K}}_w}(t,\xi,y)-p^{\widehat{{\mathfrak{K}}}_{w'}}(t,\xi,y)\right|
\leq c (|w-w'|^{\beta}\land 1) \rr_t(y-\xi)\,.
$$
Therefore,
\begin{align*}
&\left| p^{\mathfrak{K}_w}(t,x,y)- p^{\mathfrak{K}_w}(t,x',y)
-\left(p^{\mathfrak{K}_{w'}}(t,x,y)-p^{\mathfrak{K}_{w'}}(t,x',y)\right)\right|\\
& \qquad \leq  c \int_{\Rd} \left( \frac{|x-x'|}{h^{-1}(1/t)}\right) \rr_t(\xi-x)  (|w-w'|^{\beta}\land 1) \rr_t(y-\xi)\, d\xi\,.
\end{align*}
Now, by \cite[Corollary~5.14 and Lemma~5.6]{GS-2018} we get
\begin{align*}
\int_{\Rd} \rr_t(y-\xi) \rr_t(\xi-x)\, d\xi
\leq c \int_{\Rd} \rr_{2t}(y-x) \left(\rr_t(y-\xi) +\rr_t(\xi-x) \right)d\xi \leq c \rr_{2t}(y-x)\,.
\end{align*}
Finally, as in Remark~\ref{rem-r:monot_h}, we get by the monotonicity that $[h^{-1}(1/(2t))]^{-d}\leq [h^{-1}(1/t)]^{-d}$, which further gives
$\rr_{2t}(y-x)\leq 2\rr_{t}(y-x)$. This completes the proof.
\qed

\begin{lemma}\label{lem-r:cancel_01}
Let $\beta_1\in [0,\beta]\cap [0,\lah)$.
For every $T>0$ there exists a constant
$c=c(d,T,\parame,\beta_1)$
such that for all $t\in(0,T]$, $x,x' \in \Rd$,
\begin{align*}
\left| \int_{\Rd}\left( p^{\mathfrak{K}_y}(t,x,y)- p^{\mathfrak{K}_y}(t,x',y)
\right) dy\right| 
\leq c \left[ h^{-1}(1/t)\right]^{\beta_1} \left( \frac{|x-x'|}{h^{-1}(1/t)} \land 1\right).
\end{align*}
\end{lemma}
\pf
Let ${\rm I}$ be the left hand side of the inequality.
Since $\int_{\Rd}p^{\mathfrak{K}_x}(t,x,y)dy=1$, by \cite[Theorem~2.11]{GS-2018}, \cite[Proposition~5.3]{S-2018} and  \cite[Lemma~5.17(a)]{GS-2018} we have
\begin{align*}
\left| \int_{\Rd} p^{\mathfrak{K}_y}(t,x,y)\,dy - 1 \right|
&=\left| \int_{\Rd} \left( p^{\mathfrak{K}_y}(t,x,y)- p^{\mathfrak{K}_x}(t,x,y)\right)dy\right|\\
&\leq c \int_{\Rd}   (|y-x|^{\beta_1}\land 1) \rr_{t}(y-x)\, dy
\leq c \left[ h^{-1}(1/t)\right]^{\beta_1}.
\end{align*}
Now, for $|x-x'|\geq h^{-1}(1/t)$ we add and subtract 1 as follows:
\begin{align*}
{\rm I}&=\left|\left(\int_{\Rd} p^{\mathfrak{K}_y}(t,x,y)\,dy - 1 \right)
-\left(\int_{\Rd} p^{\mathfrak{K}_y}(t,x,y)\,dy - 1\right)
\right|
\leq c \left[ h^{-1}(1/t)\right]^{\beta_1}.
\end{align*}
For $|x-x'|\leq h^{-1}(1/t)$
we subtract zero, and use Lemma~\ref{lem-r:for_cancel_0}
and \cite[Lemma~5.17(a)]{GS-2018} to get
\begin{align*}
{\rm I}&=\left| \int_{\Rd}\left[
p^{\mathfrak{K}_y}(t,x,y)-p^{\mathfrak{K}_y}(t,x',y)
-\left( p^{\mathfrak{K}_x}(t,x,y)-p^{\mathfrak{K}_x}(t,x',y) \right)
\right] dy
\right|\\
&\leq c \int_{\Rd} 
\left( \frac{|x-x'|}{h^{-1}(1/t)}\right) (|y-x|^{\beta_1}\land 1) \rr_{t}(y-x)\,
  dy 
  \leq c \left[ h^{-1}(1/t)\right]^{\beta_1} \left( \frac{|x-x'|}{h^{-1}(1/t)}\right).
\end{align*}
\qed

\begin{lemma}\label{lem-r:cancel_02}
Let $r_0 \in [0,1]\cap [0,\lah+\beta\land \lah)$.
For every $T>0$ there exists a constant
$c=c(d,T,\parame,r_0)$
 such that for all $t\in (0,T]$,
$x,x',y\in\Rd$ and $r\in[0,r_0]$,
\begin{align*}
 \int_{t/2}^t \left| \int_{\Rd} \left(p^{\mathfrak{K}_z}(t-s,x,z)- p^{\mathfrak{K}_z}(t-s,x',z) \right)  q(s,z,y) \,dz\right| ds & \\
 \leq 
c \left(|x-x'|^r \land 1 \right) \left[ h^{-1}(1/t)\right]^{-1-r}&
\big( \rr_t(y-x)+ \rr_t(y-x') \big).
\end{align*}
\end{lemma}
\pf
We denote
\begin{align*}
{\rm I} &:= \left| \int_{\Rd} \left(p^{\mathfrak{K}_z}(t-s,x,z)- p^{\mathfrak{K}_z}(t-s,x',z) \right)  q(s,z,y) \,dz\right| \\
& \,\, \leq \int_{\Rd} \left| p^{\mathfrak{K}_z}(t-s,x,z)-p^{\mathfrak{K}_z}(t-s,x',z) \right| \big| q(s,z,y)-q(s,x,y) \big| \,dz \\
&\, \qquad \qquad +\, 
\left|\int_{\Rd} \left( p^{\mathfrak{K}_z}(t-s,x,z)-p^{\mathfrak{K}_z}(t-s,x',z) \right) dz \right| \left| q(s,x,y)\right|.
\end{align*}
We start by investigating
$$
{\rm I}_0:= \int_{\Rd} \left| p^{\mathfrak{K}_z}(t-s,x,z)-p^{\mathfrak{K}_z}(t-s,x',z) \right| \big| q(s,z,y)-q(s,x,y) \big| \,dz\,.
$$
By \eqref{ineq-r:est_diff_1} 
 and $|q(s,z,y)-q(s,x,y)|\leq |q(s,z,y)-q(s,x',y)|+|q(s,x',y)-q(s,x,y)|$
we get
\begin{align*}
{\rm I}_0&\leq  c \int_{\Rd} \left(\frac{|x-x'|}{h^{-1}(1/(t-s))}\land 1 \right)  \rr_{t-s}(x-z) \big| q(s,z,y)-q(s,x,y) \big| dz\\
&\quad + c \int_{\Rd} \left(\frac{|x-x'|}{h^{-1}(1/(t-s))}\land 1 \right)  \rr_{t-s}(x'-z) \big| q(s,z,y)-q(s,x',y) \big| dz\\
&\quad + c \int_{\Rd} \left(\frac{|x-x'|}{h^{-1}(1/(t-s))}\land 1 \right)  \rr_{t-s}(x'-z) \big| q(s,x,y)-q(s,x',y) \big| dz\,.
\end{align*}
Define
\begin{align*}
{\rm I}_1 
&=\int_{\Rd} (t-s)\err{0}{\beta_1-\gamma}(t-s, x-z)
\left\{\big(\err{\gamma}{0}+\err{\gamma-\beta_1}{\beta_1}\big)(s,z-y)+\big(\err{\gamma}{0}+\err{\gamma-\beta_1}{\beta_1}\big)(s,x-y)\right\}dz\,,\\
{\rm I}_2 
&=\int_{\Rd} (t-s)\err{0}{\beta_1-\gamma}(t-s, x'-z)
\left\{\big(\err{\gamma}{0}+\err{\gamma-\beta_1}{\beta_1}\big)(s,z-y)+\big(\err{\gamma}{0}+\err{\gamma-\beta_1}{\beta_1}\big)(s,x'-y)\right\}dz\,,\\
{\rm I}_3
&=(|x-x'|^{\beta_1-\gamma}\land 1)\int_{\Rd} (t-s)\err{0}{0}(t-s, x'-z)\\
& \hspace{0.25\linewidth} \times \left\{\big(\err{\gamma}{0}+\err{\gamma-\beta_1}{\beta_1}\big)(s,x-y)+
\big(\err{\gamma}{0}+\err{\gamma-\beta_1}{\beta_1}\big)(s,x'-y)\right\}dz\,.
\end{align*}
Now, let $\beta_1\in (0,\beta]\cap (0,\lah)$ be such that $\lah+\beta_1-r_0>0$
and
fix $0<\gamma_1<(\lah+\beta_1-r_0)\land \beta_1$.
By
Lemma~\ref{lem-r:q_reg_impr}
there is a constant
$c=c(d,T,\parame,\beta_1,\gamma_1)$
 such that for all $\gamma \in [\gamma_1,\beta_1]$,
$$
{\rm I}_0\leq c \left(\frac{|x-x'|}{h^{-1}(1/(t-s))}\land 1 \right) \big({\rm I}_1+{\rm I}_2 +{\rm I}_3 \big).
$$
In what follows we frequently replace $s\in (t/2,t)$ with $t$ due to
the comparability of $h^{-1}(1/s)$ and $h^{-1}(1/t)$. 
More precisely, 
using monotonicity as in Remark~\ref{rem-r:monot_h},
and additionally by \eqref{eq:intro:wlsc}, \cite[Lemma~5.3]{GS-2018} and Remark~\ref{rem-r:stretch} for the lower bound,
$$
\frac{h^{-1}(1/t)}{(2C_h [h^{-1}(1/T)\vee 1]^2)^{1/\lah}}  \leq h^{-1}(2/t) \leq h^{-1}(1/s) \leq h^{-1}(1/t)\,.
$$
Next, by \cite[Lemma~5.17(b)]{GS-2018}
with $\beta_0=\beta_1$, $m_1=n_1=n_2=\beta_1-\gamma$, $m_2=0$, we have
\begin{align*}
&\int_{\Rd} (t-s)\err{0}{\beta_1-\gamma}(t-s,x-z)
\err{\gamma}{0}(s,z-y)\,dz\\
&\leq c (t-s)\left[h^{-1}(1/s)\right]^{\gamma}
\left[\left((t-s)^{-1} \left[h^{-1}(1/(t-s))\right]^{\beta_1-\gamma}+s^{-1}\left[h^{-1}(1/s)\right]^{\beta_1-\gamma}\right) \err{0}{0}(t, x-y) \right.\\
&\hspace{0.25\linewidth} \left. + (t-s)^{-1}\left[h^{-1}(1/(t-s))\right]^{\beta_1-\gamma}\err{0}{0}(t, x-y) +s^{-1} \err{0}{\beta_1-\gamma}(t, x-y)\right]\\
&\leq c \left[h^{-1}(1/(t-s))\right]^{\beta_1-\gamma}
\err{\gamma}{0}(t, x-y)
+c (t-s)t^{-1}\big(\err{\beta_1}{0}+\err{\gamma}{\beta_1-\gamma}\big)(t, x-y)\,.
\end{align*}
By \cite[Lemma~5.17(b)]{GS-2018}
with $\beta_0=\beta_1$, $m_1=\beta_1-\gamma$, $m_2=\beta_1$, $n_1=n_2=\beta_1$, we have
\begin{align*}
&\int_{\Rd} (t-s)\err{0}{\beta_1-\gamma}(t-s, x-z)
\err{\gamma-\beta_1}{\beta_1}(s,z-y)\,dz\\
&\leq c (t-s)\left[h^{-1}(1/s)\right]^{\gamma-\beta_1}\left[\left( (t-s)^{-1}\left[h^{-1}(1/(t-s))\right]^{\beta_1}+s^{-1} \left[h^{-1}(1/s)\right]^{\beta_1}\right)  \err{0}{0}(t, x-y) \right.\\
&\hspace{0.1\linewidth} \left. + (t-s)^{-1}\left[h^{-1}(1/(t-s))\right]^{\beta_1-\gamma}\err{0}{\beta_1}(t, x-y) +s^{-1}\left[h^{-1}(1/s)\right]^{\beta_1} \err{0}{\beta_1-\gamma}(t, x-y)\right]\\
&\leq c \left[h^{-1}(1/(t-s))\right]^{\beta_1} \err{\gamma-\beta_1}{0}(t, x-y)
+ c (t-s)t^{-1}\big(\err{\gamma}{0}+\err{\gamma}{\beta_1-\gamma}\big)(t, x-y)\\
& \quad + c \left[h^{-1}(1/(t-s))\right]^{\beta_1-\gamma} \err{\gamma-\beta_1}{\beta_1}(t, x-y)\,.
\end{align*}
This gives uniformly in $\gamma\in [\gamma_1,\beta_1]$,
\begin{align*}
{\rm I}_1
&\leq c 
\left[h^{-1}(1/(t-s))\right]^{\beta_1-\gamma}
\big(\err{\gamma}{0}+\err{\gamma-\beta_1}{\beta_1}\big)(t,x-y)\\
&\quad +c \,(t-s)t^{-1}\big(\err{\beta_1}{0}+\err{\gamma}{\beta_1-\gamma}+\err{\gamma}{0}\big)(t,x-y) \\
&\quad + c \left[h^{-1}(1/(t-s))\right]^{\beta_1} \err{\gamma-\beta_1}{0}(t,x-y)\,.
\end{align*}
We treat ${\rm I}_2$ alike.
Further, by \cite[Lemma~5.6)]{GS-2018} we have
\begin{align*}
\int_{\Rd} (t-s)\err{0}{0}(t-s, x'-z)
\big(\err{\gamma}{0}+\err{\gamma-\beta_1}{\beta_1}\big)(s,x-y)\, dz
\leq c \big(\err{\gamma}{0}+\err{\gamma-\beta_1}{\beta_1}\big)(t,x-y)
\end{align*}
and so uniformly in $\gamma\in [\gamma_1,\beta_1]$,
\begin{align*}
{\rm I}_3  \leq c \,(|x-x'|^{\beta_1-\gamma}\land 1)
\left\{\big(\err{\gamma}{0}+\err{\gamma-\beta_1}{\beta_1}\big)(t,x-y)+
\big(\err{\gamma}{0}+\err{\gamma-\beta_1}{\beta_1}\big)(t,x'-y)\right\}.
\end{align*}
For each $r\in[0,r_0]$ we take $\gamma=\gamma_1\vee (\beta_1-r)$.
Using the inequalities in Remark~\ref{rem-r:monot_h}
we have
\begin{align*}
\left(\frac{|x-x'|}{h^{-1}(1/(t-s))}\land 1 \right)
&\leq \left(|x-x'|^r \land [h^{-1}(1/T)]^r \right)  \left[h^{-1}(1/(t-s))\right]^{-r} \\
&\leq [h^{-1}(1/T)\vee 1]^{r_0} \,(|x-x'|^r \land 1) \left[h^{-1}(1/(t-s))\right]^{-r}.
\end{align*}
Note that $\beta_1-\gamma-r\leq 0$ and, by considering $r\leq \beta_1-\gamma_1$ and $r>\beta_1-\gamma_1$, we obtain
$$
\frac{\beta_1-\gamma-r}{\lah}+1\geq 
\min \{1, \frac{\beta_1-\gamma_1-r_0}{\lah}+1\}>0\,.
$$
We also have $(-r/\lah) +2\geq (-r_0/\lah) +2>0$
and $((\beta_1-r)/2)\land ((\beta_1-r)/\lah)+1 \geq ((\beta_1-r_0)/2)\land ((\beta_1-r_0)/\lah)+1>0$. Therefore, \cite[Lemma~5.15]{GS-2018} and
the monotonicity of the Beta function provide, 
uniformly for all $r\in [0,r_0]$,
\begin{align*}
\int_{t/2}^t \left[h^{-1}(1/(t-s))\right]^{\beta_1-\gamma-r}ds &\leq c t \left[h^{-1}(1/t)\right]^{\beta_1-\gamma-r}\,,
\end{align*}
\begin{align*}
\int_{t/2}^t (t-s)\left[h^{-1}(1/(t-s))\right]^{-r}ds
&\leq c t^2 \left[h^{-1}(1/t)\right]^{-r}\,,
\end{align*}
\begin{align*}
\int_{t/2}^t \left[h^{-1}(1/(t-s))\right]^{\beta_1-r}ds
\leq c t \left[h^{-1}(1/t)\right]^{\beta_1-r}\,.
\end{align*}
Thus
\begin{align*}
\int_{t/2}^t \left(\frac{|x-x'|}{h^{-1}(1/(t-s))}\land 1 \right) {\rm I}_1\,ds
\leq c \,(|x-x'|^r \land 1) \,t \left[ h^{-1}(1/t)\right]^{-r}
 \big(\err{\beta_1}{0}+\err{0}{\beta_1}+\err{\gamma}{0}\big)(t,x-y)\,.
\end{align*}
We deal with 
the part containing $\rm{I}_2$ in the same way.
Similarly,
\begin{align*}
\left(\frac{|x-x'|}{h^{-1}(1/(t-s))}\land 1 \right)
\leq [h^{-1}(1/T)\vee 1]^{r_0} \,(|x-x'|^{r-(\beta_1-\gamma)} \land 1) \left[h^{-1}(1/(t-s))\right]^{\beta_1-\gamma-r}
\end{align*}
and
\begin{align*}
\int_{t/2}^t \left(\frac{|x-x'|}{h^{-1}(1/(t-s))}\land 1 \right) {\rm I}_3\,ds 
\leq c \int_{t/2}^t (|x-x'|^{r-(\beta_1-\gamma)}\land 1)\left[ h^{-1}(1/(t-s))\right]^{\beta_1-\gamma-r}\, {\rm I}_3\,ds \\
\leq c\, (|x-x'|^r \land 1 )\,t \left[ h^{-1}(1/t)\right]^{-r} \left\{\big(\err{\beta_1}{0}+\err{0}{\beta_1}\big)(t,x-y)+
\big(\err{\beta_1}{0}+\err{0}{\beta_1}\big)(t,x'-y)\right\}.
\end{align*}
To sum up, we have, uniformly for all $r\in [0,r_0]$,
\begin{align*}
\int_{t/2}^t {\rm I}_0\,ds
\leq c \left(|x-x'|^r \land 1 \right) \left[ h^{-1}(1/t)\right]^{-r}
\big( \rr_t(y-x)+ \rr_t(y-x') \big).
\end{align*}
Now, since $|q(s,x,y)|\leq c \big(\err{0}{\beta_1}+\err{\beta_1}{0}\big)(s,x-y)
\leq c \big(\err{0}{\beta_1}+\err{\beta_1}{0}\big)(t,x-y)$ (see \cite[(37)]{GS-2018}, \cite[(53)]{S-2018}), together with Lemma~\ref{lem-r:cancel_01} we get,
uniformly for all $r\in [0,r_0]$,
\begin{align*}
&\int_{t/2}^t\left|\int_{\Rd} \left( p^{\mathfrak{K}_z}(t-s,x,z)-p^{\mathfrak{K}_z}(t-s,x',z) \right) dz \right| \left| q(s,x,y)\right| ds\\
& \leq c\, (|x-x'|^{r}\land 1) \int_{t/2}^t \left[ h^{-1}(1/(t-s))\right]^{\beta_1-r}ds \, \big(\err{0}{\beta_1}+\err{\beta_1}{0}\big)(t,x-y) \\
&\leq c\, (|x-x'|^{r}\land 1) \left[h^{-1}(1/t)\right]^{-r} \big( \rr_t(y-x)+ \rr_t(y-x') \big)\,.
\end{align*}
Finally, since
\begin{align*}
\int_{t/2}^t {\rm I}\, ds
\leq \int_{t/2}^t {\rm I}_0\, ds 
+\int_{t/2}^t\left|\int_{\Rd} \left( p^{\mathfrak{K}_z}(t-s,x,z)-p^{\mathfrak{K}_z}(t-s,x',z) \right) dz \right| \left| q(s,x,y)\right| ds\,,
\end{align*}
the proof is complete.
\qed

\begin{proposition}\label{prop-r:key_0}
Let 
$r_0 \in [0,1]\cap [0,\lah+\beta\land \lah)$.
For every $T>0$ there exists a constant
$c=c(d,T,\parame,r_0)$ such that for all $t\in (0,T]$,
$x,x',y\in\Rd$ and $r\in [0,r_0]$,
\begin{align*}
\left| \phi_y(t,x)-\phi_y(t,x')\right| \leq
c \left(|x-x'|^r \land 1 \right) \left[ h^{-1}(1/t)\right]^{-r}
\big( \rr_t(y-x)+ \rr_t(y-x') \big)\,.
\end{align*}
\end{proposition}
\pf
For $s\in (0,t/2]$ we use
Corollary~\ref{cor-r:est_diff_grad}
and 
\cite[(37)]{GS-2018}, \cite[(53)]{S-2018}
to get, for all $r\in [0,1]$,
\begin{align*}
\left|\phi_y (t,x,s)-\phi_y(t,x',s) \right|
\leq c \,
(|x-x'|^{r}\land 1) \left[h^{-1}(1/(t-s))\right]^{-r}
\int_{\Rd} (t-s)
\big( \err{0}{0}(t-s, x-z)&\\
+\,\, \err{0}{0}(t-s, x'-z)\big)
\big(\err{0}{\beta_1}+\err{\beta_1}{0}\big)(s,z-y)&\,dz.
\end{align*}
Here $\beta_1\in (0,\beta]\cap (0,\lah)$ is fixed.
Since 
\cite[Lemma~5.17(b)]{GS-2018} and Remark~\ref{rem-r:monot_h} give
\begin{align*}
&\int_{\Rd} (t-s)
 \err{0}{0}(t-s, x-z)
\big(\err{0}{\beta_1}+\err{\beta_1}{0}\big)(s,z-y)\,dz\\
&\leq c \left(\left[h^{-1}(1/(t-s))\right]^{\beta_1}+\left[h^{-1}(1/s)\right]^{\beta_1}
+(t-s) s^{-1}\left[h^{-1}(1/s)\right]^{\beta_1} \right) 
\err{0}{0}(t,x-y) + \err{0}{\beta_1}(t,x-y)\\
&\leq 
c \left(\left[h^{-1}(1/t)\right]^{\beta_1}
+ t s^{-1}\left[h^{-1}(1/s)\right]^{\beta_1} \right) 
\err{0}{0}(t,x-y) + \err{0}{\beta_1}(t,x-y)\,,
\end{align*}
and \cite[Lemma~5.3]{GS-2018} gives $\left[h^{-1}(1/(t-s))\right]^{-r}\leq c \left[h^{-1}(1/t)\right]^{-r}$,
 by \cite[Lemma~5.15]{GS-2018} we get
\begin{align*}
\int_{0}^{t/2} &\left|\phi_y (t,x,s)-\phi_y(t,x',s) \right| ds\\ 
&   \leq 
c \,(|x-x'|^r \land 1) \left[ h^{-1}(1/t) \right]^{-r}
   t \left( \big( \err{\beta_1}{0}+ \err{0}{\beta_1}\big)(t,x-y)+ \big( \err{\beta_1}{0}+ \err{0}{\beta_1}\big)(t,x'-y)\right).
\end{align*}
For the remaining part of the integral
with integration in $s$ over $(t/2,t)$ we apply
Lemma~\ref{lem-r:cancel_02}.

\qed

\noindent
{\it Proof of Theorem~\ref{thm-r:1}.}
The result follows from
\eqref{e:p-kappa},
Corollary~\ref{cor-r:est_diff_grad} and
Proposition~\ref{prop-r:key_0}.
\qed

\section{Regularity - part II}

In this section  we assume that $\lah+\beta\land \lah>1$.
This condition necessitates $\lah>1/2$. 
The proofs here differ from those in Section~\ref{sec:part1}; see  Lemma~\ref{lem-r:cancel_12}.

\begin{lemma}\label{lem-r:for_cancel_1}
For every $T>0$ there exists a constant
$c=c(d,T,\parame)$
such that for all $t\in(0,T]$ and $x,x',y,w,w' \in \Rd$
satisfying $|x-x'|\leq h^{-1}(1/t)$,
\begin{align*}
\left|\nabla_x p^{\mathfrak{K}_w}(t,x,y)-\nabla_{x'} p^{\mathfrak{K}_w}(t,x',y)
-\left(\nabla_x p^{\mathfrak{K}_{w'}}(t,x,y)-\nabla_{x'} p^{\mathfrak{K}_{w'}}(t,x',y)\right)\right|&\\
\leq 
c \left( \frac{|x-x'|}{h^{-1}(1/t)}\right) (|w-w'|^{\beta}\land 1)\left[ h^{-1}(1/t)\right]^{-1} &\, \rr_{t}(y-x)\,.
\end{align*}
\end{lemma}
\pf
Let $w_0\in\Rd$  be fixed. Define $\mathfrak{K}_0(z)=(\kappa_0/(2\kappa_1)) \kappa(w_0,z)$ and $\widehat{\mathfrak{K}}_w (z)=\mathfrak{K}_w(z)- \mathfrak{K}_0(z)$. 
By \eqref{eq-r:przez_k_0-impr}, \eqref{ineq-r:est_diff_grad_1}
and
\cite[Theorem~2.11]{GS-2018}, \cite[Proposition~5.3]{S-2018},
for $|x-x'|\leq h^{-1}(1/t)$ we get
\begin{align*}
&\left|\nabla_x p^{\mathfrak{K}_w}(t,x,y)-\nabla_{x'} p^{\mathfrak{K}_w}(t,x',y)
-\left(\nabla_x p^{\mathfrak{K}_{w'}}(t,x,y)-\nabla_{x'} p^{\mathfrak{K}_{w'}}(t,x',y)\right)\right|\\
&\qquad= \left| \int_{\Rd} \left(  \nabla_x p^{\mathfrak{K}_0}(t, x,\xi)-\nabla_{x'} p^{\mathfrak{K}_0}(t, x',\xi) \right)  \left( p^{\widehat{\mathfrak{K}}_w}(t,\xi,y)-p^{\widehat{{\mathfrak{K}}}_{w'}}(t,\xi,y)\right) d\xi \right|\\
& \qquad \leq  c \int_{\Rd} \left[ h^{-1}(1/t)\right]^{-1}\left( \frac{|x-x'|}{h^{-1}(1/t)}\right) \rr_t(\xi-x)  (|w-w'|^{\beta}\land 1) \rr_t(y-\xi) d\xi\\
& \qquad\leq c \left[ h^{-1}(1/t)\right]^{-1}\left( \frac{|x-x'|}{h^{-1}(1/t)}\right) (|w-w'|^{\beta}\land 1) \rr_{2t}(y-x)\,.
\end{align*}
\cite[Corollary~5.14 and Lemma~5.6]{GS-2018}
have also been used in the last inequality.
It remains to apply $\rr_{2t}(y-x)\leq 2\rr_{t}(y-x)$; see Remark~\ref{rem-r:monot_h}.
\qed

\begin{lemma}\label{lem-r:cancel_11}
Let 
$\beta_1\in [0,\beta]\cap [0,\lah)$.
For every $T>0$ there exists a constant
$c=c(d,T,\parame,\beta_1)$
such that for all $t\in(0,T]$, $x,x' \in \Rd$,
\begin{align*}
\left| \int_{\Rd}\left(
\nabla_x p^{\mathfrak{K}_y}(t,x,y)-\nabla_{x'} p^{\mathfrak{K}_y}(t,x',y)
\right) dy \right|
\leq c \left[ h^{-1}(1/t)\right]^{-1+\beta_1} \left( \frac{|x-x'|}{h^{-1}(1/t)}\land 1\right).
\end{align*}
\end{lemma}
\pf
Let ${\rm I}$ be the left hand side of the inequality.
For $|x-x'|\geq h^{-1}(1/t)$
we conclude from \cite[(29)]{GS-2018} and \cite[Lemma~5.5]{S-2018} that
$$
{\rm I}\leq c \left[ h^{-1}(1/t)\right]^{-1+\beta_1}\,.
$$
For $|x-x'|\leq h^{-1}(1/t)$
we subtract zero, and use Lemma~\ref{lem-r:for_cancel_1}
and \cite[Lemma~5.17(a)]{GS-2018} to get
\begin{align*}
&\left| \int_{\Rd}\left[
\nabla_x p^{\mathfrak{K}_y}(t,x,y)-\nabla_{x'} p^{\mathfrak{K}_y}(t,x',y)
-\left( \nabla p^{\mathfrak{K}_x}(t,\cdot,y)(x)-\nabla_{x'} p^{\mathfrak{K}_x}(t,x',y) \right)
\right] dy
\right|\\
&\leq  c\int_{\Rd} 
\left[ h^{-1}(1/t)\right]^{-1}\left( \frac{|x-x'|}{h^{-1}(1/t)}\right) (|y-x|^{\beta_1}\land 1) \rr_{t}(y-x)\,
  dy \\
&\leq c \left[ h^{-1}(1/t)\right]^{-1+\beta_1} \left( \frac{|x-x'|}{h^{-1}(1/t)}\right)\,.
\end{align*}
\qed

\begin{lemma}\label{lem-r:cancel_12}
Let  $r_0\in [0,1]\cap [0,\lah+\beta\land \lah-1)$.
For every $T>0$ there exists a constant
$c=c(d,T,\parame, r_0)$
such that for all $t\in (0, T]$, $x,x' \in \Rd$ and $r\in [0,r_0]$,
\begin{align*}
 \int_{t/2}^t \left| \int_{\Rd} \left(\nabla_x p^{\mathfrak{K}_z}(t-s,x,z)-\nabla_{x'} p^{\mathfrak{K}_z}(t-s,x',z) \right)  q(s,z,y) \,dz\right| ds & \\
 \leq 
c \left(|x-x'|^r \land 1 \right) \left[ h^{-1}(1/t)\right]^{-1-r}&
\big( \rr_t(y-x)+ \rr_t(y-x') \big).
\end{align*}
\end{lemma}
\pf
Denote
\begin{align*}
{\rm I}:= \left| \int_{\Rd} \left(\nabla_x p^{\mathfrak{K}_z}(t-s,x,z)-\nabla_{x'} p^{\mathfrak{K}_z}(t-s,x',z) \right)  q(s,z,y) \,dz\right|.
\end{align*}
Let $\beta_1\in (0,\beta]\cap (0,\lah)$ 
be such that $\lah+\beta_1-1-r_0>0$.
Fix $0<\gamma < (\lah+\beta_1-1-r_0)\land \beta_1$.
We fist show that
\begin{align}\label{eq-r:raw}
{\rm I}\leq c  \left( \frac{|x-x'|}{h^{-1}(1/(t-s))}\land 1\right)
\Big( V (t,x-y;s)+ V(t,x'-y;s) \Big).
\end{align}
where
\begin{align*}
V(t,x-y;s):=\left[ h^{-1}(1/(t-s))\right]^{-1} \Big\{& \left[h^{-1}(1/(t-s))\right]^{\beta_1-\gamma} \big(\err{\gamma}{0}+\err{\gamma-\beta_1}{\beta_1}\big)(t,x-y)\\
& +(t-s)t^{-1} \big(\err{\beta_1}{0}+\err{\gamma}{\beta_1-\gamma} + \err{\gamma}{0}\big)(t,x-y)\\
& +\left[h^{-1}(1/(t-s))\right]^{\beta_1} \big(\err{\beta_1}{0}+\err{0}{\beta_1}+\err{\gamma-\beta_1}{0}\big)(t,x-y) \Big\}.
\end{align*}
In what follows we frequently replace $s\in (t/2,t)$ with $t$ due to
the comparability of $h^{-1}(1/s)$ with $h^{-1}(1/t)$; see
\cite[Lemmas 5.1 and~5.3]{GS-2018}
and Remark~\ref{rem-r:stretch}.
For $|x-x'|\geq  h^{-1}(1/(t-s))$ we have,
by   \cite[Proposition~2.1, (38), (29) and~(37)]{GS-2018} and
\cite[Proposition~4.1, (54), Lemma~5.5 and~(53)]{S-2018},
\begin{align*}
&\left| \int_{\Rd}\nabla_x p^{\mathfrak{K}_z}(t-s,x,z)  q(s,z,y) \,dz\right|\\
&\leq  \int_{\Rd} \left| \nabla_x p^{\mathfrak{K}_z}(t-s,x,z) \right| |  q(s,z,y)-q(s,x,y)| \,dz
 + \left| \int_{\Rd}\nabla_x p^{\mathfrak{K}_z}(t-s,x,z)\,dz\right|
|q(s,x,y)|\\
&\leq c \left[ h^{-1}(1/(t-s))\right]^{-1} \bigg\{ \int_{\Rd}(t-s)\err{0}{\beta_1-\gamma}(t-s, x-z)
\big(\err{\gamma}{0}+\err{\gamma-\beta_1}{\beta_1}\big)(s,z-y)\,dz \\
&\hspace{0.26\linewidth}+\int_{\Rd}(t-s)\err{0}{\beta_1-\gamma}(t-s, x-z)\,dz \,\big(\err{\gamma}{0}+\err{\gamma-\beta_1}{\beta_1}\big)(t,x-y)\\
&\hspace{0.26\linewidth} + \left[ h^{-1}(1/(t-s))\right]^{\beta_1} \big(\err{\beta_1}{0}+\err{0}{\beta_1}\big)(t,x-y)\bigg\} =: {\rm R}.
\end{align*}
Now, \eqref{eq-r:raw} follows in this case from \cite[Lemma~5.18(a) and (b)]{GS-2018} (once with $n_1=n_2=\beta_1$).
For $|x-x'|\leq h^{-1}(1/(t-s))$, 
by \eqref{ineq-r:est_diff_grad_1},
Lemma~\ref{lem-r:cancel_11} and 
\cite[(38), (37)]{GS-2018}, \cite[(54), (53)]{S-2018}
we have
\begin{align*}
{\rm I}&\leq \int_{\Rd} \left| \nabla_x p^{\mathfrak{K}_z}(t-s,x,z)-\nabla_{x'} p^{\mathfrak{K}_z}(t-s,x',z) \right| | q(s,z,y)-q(s,x,y)| \,dz\\
&\quad + \left| \int_{\Rd} \left(\nabla_x p^{\mathfrak{K}_z}(t-s,x,z)-\nabla_{x'} p^{\mathfrak{K}_z}(t-s,x',z)\right) \,dz \right|  |q(s,x,y)|
\leq c
\left( \frac{|x-x'|}{h^{-1}(1/(t-s))}\right) {\rm R} \,.
\end{align*}
Here again \eqref{eq-r:raw} follows from \cite[Lemma~5.18(a) and (b)]{GS-2018}.
Finally, 
since by our assumptions $(\beta_1-\gamma-1-r)/\lah+1\geq (\beta_1-\gamma-1-r_0)/\lah+1>0$
and $(-1-r)/\lah +2\geq (-1-r_0)/\lah +2>0$, inequality \eqref{eq-r:raw} and \cite[Lemma~5.15]{GS-2018}
with the monotonicity of Beta function give, uniformly for all $r\in [0,r_0]$,
\begin{align*}
&\int_{t/2}^t {\rm I} \,ds\leq c \,(|x-x'|^r \land 1) \int_{t/2}^t
\left[h^{-1}(1/(t-s))\right]^{-r} \Big( V(t,x-y;s)+ V(t,x'-y;s) \Big)\,ds\\
&\leq c \,(|x-x'|^r \land 1)\, t \left[h^{-1}(1/t)\right]^{-1-r}\Big\{\big(\err{\beta_1}{0}+\err{0}{\beta_1}+\err{\gamma}{0}\big)(t,x-y)+
\big(\err{\beta_1}{0}+\err{0}{\beta_1}+\err{\gamma}{0}\big)(t,x'-y)\Big\}.
\end{align*}
This ends the proof (see Remark~\ref{rem-r:monot_h}).
\qed

\begin{proposition}\label{prop-r:key_1}
Let  $r_0\in [0,1]\cap [0,\lah+\beta\land \lah-1)$.
For every $T>0$ there exists a constant
$c=c(d,T,\parame, r_0)$ such that for all $t\in (0,T]$,
$x,x',y\in\Rd$ and $r\in [0,r_0]$,
\begin{align*}
\left|\nabla_x \phi_y(t,x)-\nabla_{x'}\phi_y(t,x')\right| \leq
c \left(|x-x'|^r \land 1 \right) \left[ h^{-1}(1/t)\right]^{-1-r}
\big( \rr_t(y-x)+ \rr_t(y-x') \big)\,.
\end{align*}
\end{proposition}
\pf
Applying \cite[(43) and~(45)]{GS-2018} and \cite[(59) and~(56)]{S-2018} we have
\begin{align*}
&\nabla_x \phi_y(t,x)-\nabla_{x'}\phi_y(t,x')
=\int_0^t \big( \nabla_x \phi_y (t,x,s)-\nabla_{x'} \phi_y(t,x',s) \big)\,ds\\
&\quad = \int_{0}^t \int_{\Rd} \left(\nabla_x p^{\mathfrak{K}_z}(t-s,x,z)-\nabla_{x'} p^{\mathfrak{K}_z}(t-s,x',z) \right)  q(s,z,y) \,dz ds\,.
\end{align*}
For $s\in (0,t/2]$ we find by 
Corollary~\ref{cor-r:est_diff_grad}
and \cite[(37)]{GS-2018}, \cite[(53)]{S-2018}
that for all $r\in [0,1]$,
\begin{align*}
&\left| \nabla_x \phi_y (t,x,s)-\nabla_{x'} \phi_y(t,x',s) \right|
\leq c \,
(|x-x'|^{r}\land 1) \left[h^{-1}(1/(t-s))\right]^{-1-r}\\
&\qquad \times \int_{\Rd} (t-s)
\big( \err{0}{0}(t-s, x-z)
+\,\, \err{0}{0}(t-s, x'-z)\big)
\big(\err{0}{\beta_1}+\err{\beta_1}{0}\big)(s,z-y)\,dz\,,
\end{align*}
where $\beta_1\in (0,\beta]\cap (0,\lah)$ is fixed.
The rest of this part of the proof is the same as that of
Proposition~\ref{prop-r:key_0}, and relies on \cite[Lemmas~5.17(b), 5.3 and~5.15]{GS-2018}, integration in $s\in (0,t/2]$ and Remark~\ref{rem-r:monot_h}.
For integration in $s$ over $(t/2,t)$ we apply
Lemma~\ref{lem-r:cancel_12}.
\qed

\noindent
{\it Proof of Theorem~\ref{thm-r:2a}.}
The result follows from
\eqref{e:p-kappa},
Corollary~\ref{cor-r:est_diff_grad} and
Proposition~\ref{prop-r:key_1}.
\qed

\section{Regularity - part III}

In this section  we assume that  $\lah+\beta\land \lah>2$.
Note that this may only hold if $\lah>1$, which in turn puts us into case $\Pa$ or $\Pc$.
We first prove that the second oder derivatives  of $p^{\kappa}(t,x,y)$ in $x$ actually exist.
Such a result is missing in \cite{GS-2018}, therefore we first
need to prepare  several technical lemmas.

\begin{lemma}\label{lem-r:for_cancel_2}
For every $T>0$ there exists a constant
$c=c(d,T,\parame)$
such that for all $t\in(0,T]$ and $x,x',y,w,w' \in \Rd$,
\begin{align}\label{ineq-r:2a}
\left|\nabla_x^2\, p^{\mathfrak{K}_w}(t,x,y)
-\nabla_x^2\, p^{\mathfrak{K}_{w'}}(t,x,y)\right|
\leq 
c\, (|w-w'|^{\beta}\land 1)\left[ h^{-1}(1/t)\right]^{-2} \, \rr_{t}(y-x)\,,
\end{align}
and if $|x-x'|\leq h^{-1}(1/t)$, then
\begin{align}\label{ineq-r:2b}
\left|\nabla_x^2\, p^{\mathfrak{K}_w}(t,x,y)-\nabla_{x'}^2\, p^{\mathfrak{K}_w}(t,x',y)
-\left(\nabla_x^2\, p^{\mathfrak{K}_{w'}}(t,x,y)-\nabla_{x'}^2\, p^{\mathfrak{K}_{w'}}(t,x',y)\right)\right|&\\
\nonumber \leq 
c \left( \frac{|x-x'|}{h^{-1}(1/t)}\right) (|w-w'|^{\beta}\land 1)\left[ h^{-1}(1/t)\right]^{-2} &\, \rr_{t}(y-x)\,.
\end{align}
\end{lemma}
\pf
By \eqref{eq-r:przez_k_0-impr} 
(\eqref{ineq-r:est_diff_1} and \eqref{ineq-r:est_diff_grad_1} allow differentiating under the integral)
we have
\begin{align*}
\nabla_x^2\, p^{\mathfrak{K}_w}(t,x,y)
-\nabla_x^2\, p^{\mathfrak{K}_{w'}}(t,x,y)=
\int_{\Rd}\nabla_x^2\, p^{\mathfrak{K}_0}(t, x,\xi) \left(p^{\widehat{\mathfrak{K}}_w}(t,\xi,y)-p^{\widehat{{\mathfrak{K}}}_{w'}}(t,\xi,y)\right)d\xi\,,
\end{align*}
and
\begin{align*}
&\left|\nabla_x^2\, p^{\mathfrak{K}_w}(t,x,y)-\nabla_{x'}^2\, p^{\mathfrak{K}_w}(t,x',y)
-\left(\nabla_x^2\, p^{\mathfrak{K}_{w'}}(t,x,y)-\nabla_{x'}^2\, p^{\mathfrak{K}_{w'}}(t,x',y)\right)\right|\\
&\qquad= \left| \int_{\Rd} \left(  \nabla_x^2\, p^{\mathfrak{K}_0}(t, x,\xi)-\nabla_{x'}^2\, p^{\mathfrak{K}_0}(t, x',\xi) \right)  \left( p^{\widehat{\mathfrak{K}}_w}(t,\xi,y)-p^{\widehat{{\mathfrak{K}}}_{w'}}(t,\xi,y)\right) d\xi \right|\,.
\end{align*}
The desired inequalities  follow from
\cite[Proposition~2.1, Theorem~2.11, Corollary~5.14 and Lemma~5.6]{GS-2018} and \eqref{ineq-r:est_diff_grad_2}; cf. proof of Lemma~\ref{lem-r:for_cancel_1}.
\qed

\begin{lemma}
Let 
$\beta_1\in [0,\beta]\cap [0,\lah)$.
For every $T>0$ there exists a constant
$c=c(d,T,\parame,\beta_1)$
such that for all $t\in(0,T]$ and $x,x' \in \Rd$,
\begin{align}\label{ineq-r:2a_int}
\left|\int_{\Rd}\nabla_x^2\, p^{\mathfrak{K}_y}(t,x,y) \,dy \right|\leq c  \left[ h^{-1}(1/t)\right]^{-2+\beta_1}\,,
\end{align}
\begin{align}\label{ineq-r:2b_int}
\left| \int_{\Rd}\left(
\nabla_x^2\, p^{\mathfrak{K}_y}(t,x,y)-\nabla_{x'}^2\, p^{\mathfrak{K}_y}(t,x',y)
\right) dy \right|
\leq c \left[ h^{-1}(1/t)\right]^{-2+\beta_1} \left( \frac{|x-x'|}{h^{-1}(1/t)}\land 1\right).
\end{align}
\end{lemma}
\pf
The proof of \eqref{ineq-r:2a_int}
is like that of 
\cite[(29)]{GS-2018},
but it requires the use of \eqref{ineq-r:2a}. For the proof of \eqref{ineq-r:2b_int} we employ
\eqref{ineq-r:2a_int} if $|x-x'|\geq h^{-1}(1/t)$,
and we use 
\eqref{ineq-r:2b} if $|x-x'|\leq h^{-1}(1/t)$; cf. proof of
Lemma~\ref{lem-r:cancel_11}.
\qed

\begin{lemma}\label{lem-r:phi_sec_der_1}
For all $0<s<t$ and $x,y\in\Rd$,
\begin{align*}
\nabla_x^2\, \phi_{y}(t,x,s)=\int_{\Rd}
\nabla_x^2\, p^{\mathfrak{K}_z}(t-s,x,z)q(s,z,y)\,dz\,.
\end{align*}
\end{lemma}
\pf By
\cite[(45)]{GS-2018} we have
$
\nabla_x \phi_{y}(t,x,s)=\int_{\Rd}
\nabla_x p^{\mathfrak{K}_z}(t-s,x,z)q(s,z,y)\,dz
$.
We obtain the result by
applying this equality to the difference quotient
$(\nabla_x \phi_{y}(t,x+\varepsilon e_i,s)-\nabla_x \phi_{y}(t,x,s))/\varepsilon$
and using  the dominated convergence theorem
justified by
\eqref{ineq-r:est_diff_grad_1} and
\cite[(37), Lemma~5.17(b)]{GS-2018}.

\qed

\begin{proposition}\label{prop-r:phi_sec_der_2}
For every $T>0$ there exists a constant
$c=c(d,T,\parame)$
such that for all $t\in(0,T]$ and $x,y\in\Rd$,
\begin{align}
\nabla_x^2\, \phi_{y}(t,x)&=\int_0^t \int_{\Rd}
\nabla_x^2 p^{\mathfrak{K}_z}(t-s,x,z)q(s,z,y)\,dzds\,, \label{eq-r:sec_der_2} \\
&\nonumber\\
| \nabla_x^2 \,\phi_y(t,x) |
& \leq c \left[ h^{-1}(1/t)\right]^{-2}\rr_t (x-y)\,.\label{ineq-r:sec_der_2}
\end{align}
\end{proposition}
\pf
We choose $\beta_1\in (0,\beta]\cap (0,\lah)$ 
such that $\lah +\beta_1-2>0$.
Let 
$0<|\varepsilon| \leq h^{-1}(3/t)$
 and $\widetilde{x}=x+\varepsilon\theta e_i$.
Based on 
\cite[Theorem~7.21]{MR924157}
and Lemma~\ref{lem-r:phi_sec_der_1}
we have
\begin{align*}
&{\rm I}_0:=\left| \frac{1}{\varepsilon} \left(\frac{\partial}{\partial x_j} \phi_y(t,x+\varepsilon e_i,s)- \frac{\partial}{\partial x_j}\phi_y(t,x,s)\right)\right|
= \left| \int_0^1 \int_{\Rd} \frac{\partial^2}{\partial x_i \partial x_j} \, p^{\mathfrak{K}_z}(t-s,\widetilde{x},z) q(s,z,y)\, dzd\theta \right|.
\end{align*}
For $s\in (0,t/2]$ by \cite[Proposition~2.1, (37), Lemmas~5.17(b) and~5.3, and Proposition~5.8]{GS-2018} and Remark~\ref{rem-r:monot_h},
\begin{align*}
{\rm I}_0&\leq c \int_0^1  \int_{\Rd} (t-s) 
\err{-2}{0}(t-s, \widetilde{x}-z) 
 \big(\err{0}{\beta_1}+\err{\beta_1}{0}\big)(s,z-y)  \,dzd\theta\\
&\leq c \left[h^{-1}(1/(t-s))\right]^{-2} \left(\int_0^1 \err{0}{0}(t, \widetilde{x}-y) \,d\theta \right)
\left(1 
+\left[h^{-1}(1/s)\right]^{\beta_1}
 + (t-s) s^{-1}\left[h^{-1}(1/s)\right]^{\beta_1}\right) \\
&\leq c \left[ h^{-1}(1/t) \right]^{-2}
\err{0}{0}(t,x-y)
\left(1+ t \,s^{-1}\left[h^{-1}(1/s)\right]^{\beta_1}\right).
\end{align*}
Next, for $s\in (t/2,t)$ we take $0<\gamma<(\lah+\beta_1-2)\land \beta_1$
and by \cite[Proposition~2.1, (38), (37)]{GS-2018}
and \eqref{ineq-r:2a_int}
we obtain
\begin{align*}
{\rm I}_0& \leq \int_0^1 \int_{\Rd} \left| \frac{\partial^2}{\partial x_i \partial x_j}\, p^{\mathfrak{K}_z}(t-s,\widetilde{x},z)\right| |q(s,z,y)-q(s,\widetilde{x},y)|\, dzd\theta\\
&\quad+ \int_0^1 \left| \int_{\Rd}   \frac{\partial^2}{\partial x_i \partial x_j}\, p^{\mathfrak{K}_z}(t-s,\widetilde{x},z)\,dz\right| |q(s,\widetilde{x},y)| \,d\theta\\
&\leq 
c \int_0^1 \int_{\Rd} (t-s) \err{-2}{\beta_1-\gamma} (t-s,\widetilde{x}-z) 
\big(\err{\gamma}{0}+\err{\gamma-\beta_1}{\beta_1}\big)(s,z-y)\,
 dzd\theta\\
&\quad+c \int_0^1 \left( \int_{\Rd} (t-s) \err{-2}{\beta_1-\gamma} (t-s,\widetilde{x}-z)\,dz\right) \big(\err{\gamma}{0}+\err{\gamma-\beta_1}{\beta_1}\big)(s,\widetilde{x}-y)\, d\theta
\\
&\quad+c \left[h^{-1}(1/(t-s))\right]^{-2+\beta_1} \int_0^1 \big(\err{0}{\beta_1}+\err{\beta_1}{0}\big)(s,\widetilde{x}-y)\,d\theta
=: {\rm I}_1+{\rm I}_2+{\rm I}_3\,.
\end{align*}
By
\cite[Lemma~5.17(a) and~(b), Lemma~5.3 and Proposition~5.8]{GS-2018}
and Remark~\ref{rem-r:monot_h} we get
\begin{align*}
{\rm I}_1\leq c\left (  \left[h^{-1}(1/(t-s))\right]^{-2+\beta_1-\gamma} \err{\gamma-\beta_1}{0}(t,x-y)+ (t-s) \left[h^{-1}(1/(t-s))\right]^{-2}t^{-1}\err{0}{0}(t,x-y) \right)\!,
\end{align*}
\begin{align*}
{\rm I}_2\leq c \left[h^{-1}(1/(t-s))\right]^{-2+\beta_1-\gamma} \err{\gamma-\beta_1}{0}(t,x-y)\,,
\qquad
{\rm I}_3 \leq c \left[h^{-1}(1/(t-s))\right]^{-2+\beta_1}\err{0}{0}(t,x-y)\,.
\end{align*}
Thus ${\rm I}_0$ is bounded by a function independent of $\varepsilon$,
which is integrable in $s$ over $(0,t)$ due to \cite[Lemma~5.15]{GS-2018}
and our assumptions that guarantee $(\beta_1-\gamma-2)/\lah+1>0$
and $(-2)/\lah+2>0$.
Now, by
\cite[(43) and~(45)]{GS-2018}
we have
\begin{align*}
 \frac{1}{\varepsilon} \left(\frac{\partial}{\partial x_j} \phi_y(t,x+\varepsilon e_i)- \frac{\partial}{\partial x_j}\phi_y(t,x)\right)
=\int_0^t \frac{1}{\varepsilon} \left(\frac{\partial}{\partial x_j} \phi_y(t,x+\varepsilon e_i,s)- \frac{\partial}{\partial x_j}\phi_y(t,x,s)\right) ds\,,
\end{align*}
and we use the dominated convergence theorem
and Lemma~\ref{lem-r:phi_sec_der_1} to reach \eqref{eq-r:sec_der_2}.
The estimate \eqref{ineq-r:sec_der_2} follows from
integrating the upper bound of ${\rm I}_0$ and applying \cite[Lemma~5.15]{GS-2018}.
\qed

We now concentrate on the regularity of $\nabla_x^2\,\phi_y(t,x)$.

\begin{lemma}\label{lem-r:cancel_22}
Let  $r_0\in [0,1]\cap [0,\lah+\beta\land \lah-2)$.
For every $T>0$ there exists a constant
$c=c(d,T,\parame, r_0)$
such that for all $t\in (0, T]$, $x,x' \in \Rd$ and $r\in [0,r_0]$,
\begin{align*}
 \int_{t/2}^t \left| \int_{\Rd} \left(\nabla_x^2 p^{\mathfrak{K}_z}(t-s,x,z)-\nabla_{x'}^2 p^{\mathfrak{K}_z}(t-s,x',z) \right)  q(s,z,y) \,dz\right| ds & \\
 \leq 
c \left(|x-x'|^r \land 1 \right) \left[ h^{-1}(1/t)\right]^{-2-r}&
\big( \rr_t(y-x)+ \rr_t(y-x') \big).
\end{align*}
\end{lemma}
\pf
The proof goes by the same lines as the proof of Lemma~\ref{lem-r:cancel_12} with $1$ replaced by $2$ in the choice of $\beta_1$ and $\gamma$, $[h^{-1}(1/u)]^{-1}$ replaced by $[h^{-1}(1/u)]^{-2}$,
\eqref{ineq-r:est_diff_grad_1} by \eqref{ineq-r:est_diff_grad_2},
\cite[(29)]{GS-2018} by \eqref{ineq-r:2a_int},
and Lemma~\ref{lem-r:cancel_11} by \eqref{ineq-r:2b_int}.
Note also that by our assumptions,
 $(\beta_1-\gamma-2-r)/\lah+1\geq (\beta_1-\gamma-2-r_0)/\lah+1>0$
and $(-2-r)/\lah +2\geq (-2-r_0)/\lah +2>0$.
\qed

\begin{proposition}\label{prop-r:key_2}
Let  $r_0\in [0,1]\cap [0,\lah+\beta\land \lah-2)$.
For every $T>0$ there exists a constant
$c=c(d,T,\parame, r_0)$ such that for all $t\in (0,T]$,
$x,x',y\in\Rd$ and $r\in [0,r_0]$,
\begin{align*}
\left|\nabla_x^2 \phi_y(t,x)-\nabla_{x'}^2\phi_y(t,x')\right| \leq
c \left(|x-x'|^r \land 1 \right) \left[ h^{-1}(1/t)\right]^{-2-r}
\big( \rr_t(y-x)+ \rr_t(y-x') \big)\,.
\end{align*}
\end{proposition}
\pf
The result follows from
\eqref{eq-r:sec_der_2}, Lemma~\ref{lem-r:phi_sec_der_1}, Corollary~\ref{cor-r:est_diff_grad}, \cite[(37), Lemmas~5.17(b), 5.3 and~5.15]{GS-2018}, integration in $s\in (0,t/2]$, Remark~\ref{rem-r:monot_h}
and Lemma~\ref{lem-r:cancel_22};
cf. proof of  Proposition~\ref{prop-r:key_1}.

\qed

\noindent
{\it Proof of Theorem~\ref{thm-r:2}.}
From
\eqref{e:p-kappa}
and \eqref{eq-r:sec_der_2}
we get the second order differentiability of 
$p^{\kappa}(t,x,y)$ in $x$.
By \cite[Proposition~2.1]{GS-2018}
and \eqref{ineq-r:sec_der_2}
we obtain the upper bound.
Finally, 
Corollary~\ref{cor-r:est_diff_grad} and
Proposition~\ref{prop-r:key_2} 
give the regularity.
\qed

\section*{Acknowledgement}

The author thanks Moritz Kassmann and Aleksei Kulik
for pointing out the issue of regularity of the heat kernels constructed in \cite{GS-2018}. The author also thanks
Krzysztof Bogdan and Tomasz Grzywny for discussions and suggestions.

\small
\bibliographystyle{abbrv}

\end{document}